\documentclass[pdflatex,sn-mathphys-num]{sn-jnl}

\usepackage[T1]{fontenc}
\usepackage[utf8]{inputenc}
\usepackage{graphicx}
\usepackage{amsmath,amssymb,amsfonts,amstext,amscd}
\usepackage{mathrsfs}
\usepackage{amsthm}
\usepackage{enumerate}
\usepackage{indentfirst}
\usepackage{xcolor}
\usepackage{url}
\usepackage[all]{xy}

\newcommand{\Z}{\mathbb{Z}}
\newcommand{\Q}{\mathbb{Q}}

\newcommand{\Fs}{\dot{F}^2}

\newcommand{\dtF}{\dot{F}}

\newcommand{\twotorsion}{{}_2\mathrm{Br}(F)}

\newcommand{\Br}[1]{\mathrm{Br}(#1)}

\makeatletter
\@ifundefined{th@thmstyleone}{%
  \newtheoremstyle{thmstyleone}%
  {12pt}{12pt}{\itshape}{0pt}{\bfseries}{}{0.5em}{}%
}{}
\@ifundefined{th@thmstylethree}{%
  \newtheoremstyle{thmstylethree}%
  {12pt}{12pt}{\normalfont}{0pt}{\bfseries}{}{0.5em}{}%
}{}
\makeatother

\theoremstyle{thmstyleone}
\newtheorem{thm}{Theorem}
\newtheorem*{theorem*}{Theorem}
\newtheorem{prop}[thm]{Proposition}

\newtheorem{cor}[thm]{Corollary}

\theoremstyle{thmstylethree}
\newtheorem{defn}[thm]{Definition}

\newtheorem{rem}[thm]{Remark}
\newtheorem{ex}[thm]{Example}

\newtheorem{conjec}{Conjecture}

\begin{document}

\title[Free Pro-$2$ Product Decompositions and Hilbert's Theorem 90]{Free Pro-$2$ Product Decompositions and Hilbert's Theorem 90 for the Kaplansky Radical}

\author*[1]{\fnm{Ronie Peterson} \sur{Dario}}\email{ronie@utfpr.edu.br}

\affil*[1]{\orgdiv{Department of Mathematics},
\orgname{Federal University of Technology -- Paran\'a (UTFPR)},
\orgaddress{\street{Av. Sete de Setembro, 3165}, \city{Curitiba},
\postcode{80230-901}, \state{Paran\'a}, \country{Brazil}}}

\abstract{The H-conjecture of Kijima and Nishi predicts a Hilbert Theorem~90 type property for the Kaplansky radical with respect to quadratic field extensions. In this paper, we establish a Galois-theoretic criterion ensuring this property for fields $F$ with finitely many square classes. The criterion is formulated in terms of a free pro-$2$ product decomposition of the maximal pro-$2$ Galois group $G_F(2)$.

As an application, we prove the H-conjecture for fields whose maximal pro-$2$ Galois groups are of elementary type, in the sense of the pro-$2$ version of Efrat's Elementary Type Conjecture.

Our criterion also applies to pre-$2$-henselian fields, to the previously known positive cases established by Iwakami, Kijima, and Nishi, and to the families of fields with nontrivial Kaplansky radical constructed by Berman and Kula. It thereby provides a unified Galois-theoretic explanation for the currently known positive instances of the finite-dimensional H-conjecture.}


\maketitle

\section{Introduction}

For a field $F$ of characteristic different from $2$, the Kaplansky radical
$R(F)$ can be defined as
\begin{equation}\label{eq:radical-as-intersection}
R(F)=\bigcap_{a\in F\setminus \{0\}}D_F\langle1,a\rangle,
\end{equation}
where $D_F\langle1,a\rangle$ denotes the set of nonzero values
represented by the binary quadratic form $X^2+aY^2$ \cite{kaplansky1}. It arises naturally in the algebraic theory of quadratic forms and Witt rings; see \cite[Chapter XII, \S6]{lam1}.

Let $\Fs = \{\alpha^2\mid\alpha\in\dot F\}$ denote the squares subgroup of $\dot F = F\setminus \{0\}$. Every square is represented by $X^2+aY^2$ for every $a\in\dot F$, and therefore $\Fs\subseteq R(F)$. This inclusion may be strict. This phenomenon was first emphasized by Cordes \cite{cordes1}, who observed that fields satisfying \begin{equation}\label{eq: nontrivial definition} \dot F^2\subsetneq R(F)\end{equation}
were necessary for a complete classification of the Witt rings of
fields with eight square classes. A field satisfying
\eqref{eq: nontrivial definition} is said to have a
\emph{nontrivial Kaplansky radical}. Subsequently, Berman
\cite{berman} and Kula \cite{kula1,kula2} constructed several families
of such fields.

The existence of fields with nontrivial Kaplansky radical gives rise to
natural analogues of classical results involving square classes,
obtained by replacing the subgroup $\dot F^2$ with $R(F)$, as in the
theory of pre-$2$-henselian valuations \cite{dario3}.   A central example of this principle is the following finite-dimensional version of the H-conjecture.

\begin{conjec}[H-conjecture -- finite-dimensional case]
\label{claimh90}
Let $F$ be a field of characteristic different from $2$ such that
$\dot F/\Fs$ is finite. Let
$K=F(\sqrt a)$, where $a\in\dot F\setminus\Fs$, and let
$N=N_{K/F}:\dot K\to\dot F$ denote the norm map. Then
\[
N^{-1}(R(F))=\dot F\,R(K).
\]
\end{conjec}

Kijima and Nishi \cite{kijima1} originally proposed the H-conjecture
without the finiteness assumption and established it for
quasi-Pythagorean fields \cite{kijima2}. The terminology reflects the
fact that for the quadratic extension $K=F(\sqrt a)$ with norm map
$N=N_{K/F}$, Hilbert's Theorem~90 implies
\begin{equation}\label{eq:H90-quadratic}
N^{-1}(\dot F^2)=\dot F\dot K^2.
\end{equation}

Becher and Leep \cite[Theorem~4.8]{becher-leep} constructed a quadratic
extension $K/F$ such that $\dot F\,R(K)\subsetneq N^{-1}(R(F)),$ 
thereby disproving the original H-conjecture when no finiteness
condition is imposed. The finite-dimensional case, however, remains
open.

The aim of this paper is to show that the norm identity in
Conjecture~\ref{claimh90} follows from a suitable free pro-$2$ product
decomposition of the maximal pro-$2$ Galois group of $F$.

Let $F(2)$ denote the maximal $2$-extension of $F$ contained in a fixed
separable closure, and write
\[G_F(2)=\operatorname{Gal}(F(2)/F)\]
for the maximal pro-$2$ Galois group of $F$. Our main result is the
following.

\begin{theorem*}[Theorem~\ref{thm:H90-decomposition}] Let $F$ be a field of characteristic different from $2$. Assume that there exist intermediate fields $F_0,F_1$ of $F(2)/F$ such that 	
\[G_F(2)\cong G_{F_0}(2)\ast_2G_{F_1}(2),\]
where $G_{F_0}(2)$ is a free pro-$2$ group and $F_1$ has trivial Kaplansky radical. Let $K=F(\sqrt a)$, where $a\in\dot F\setminus\dot F^2$, and let
$N=N_{K/F}:\dot K\longrightarrow\dot F$ be the norm map. Then
\[N^{-1}(R(F))=\dot FR(K).\]
\end{theorem*}

A principal source of fields satisfying the hypotheses of
Theorem~\ref{thm:H90-decomposition} is provided by the elementary type
class of cyclotomic pairs occurring in Efrat's Elementary Type
Conjecture \cite{efrat-livro}.

In the pro-$2$ setting considered here
\cite[Section~4]{ware3}, a cyclotomic pro-$2$ pair is a pair
$(G,\theta)$, where $G$ is a pro-$2$ group and $\theta:G\longrightarrow\mathbb Z_2^\times$ 
is a continuous homomorphism. Such a pair is called realizable if it
is isomorphic to the canonical cyclotomic pro-$2$ pair $(G_F(2),\theta_F)$ 
of some field $F$; see Definition~\ref{def:par-ciclotomico}.

The elementary type class $\mathcal C$, introduced in
Definition~\ref{def:elementary-type-class}, is generated from two
families of basic realizable pairs. The first consists of the canonical
pairs $(G_F(2),\theta_F)$ for which $G_F(2)$ is a finitely generated
free pro-$2$ group. The second consists of the canonical pairs
associated with $2$-Demushkin fields, in the sense of
Definition~\ref{def:demushkin-field-and-group}. The class $\mathcal C$
is the smallest class of realizable cyclotomic pro-$2$ pairs containing
these basic pairs and closed under free pro-$2$ products and
cyclotomic semidirect products.

The pro-$2$ version of Efrat's Elementary Type Conjecture predicts
that every realizable cyclotomic pro-$2$ pair $(G,\theta)$ whose underlying pro-$2$ group $G$ is topologically finitely generated belongs to $\mathcal C$. The construction of $\mathcal C$ and the precise
statement of the conjecture are recalled in
Subsection~\ref{subsec:ETG} of
Section~\ref{sec:the-elementary-case}.

Theorem~\ref{thm:claim2-ETG2} shows that, whenever $F$ is a field
whose canonical cyclotomic pro-$2$ pair $(G_F(2),\theta_F)$ belongs to
$\mathcal C$, there exist intermediate fields $L,H$ of $F(2)/F$ such
that, in the category of cyclotomic pairs,
\[(G_F(2),\theta_F)\cong(G_L(2),\theta_L)\ast_2(G_H(2),\theta_H),\]
where $G_L(2)$ is free pro-$2$ and $H$ has trivial Kaplansky radical.
Consequently, the corresponding decomposition of the underlying
pro-$2$ groups satisfies the hypotheses of Theorem~\ref{thm:H90-decomposition}. We therefore obtain the following
consequence.

\begin{theorem*}[Theorem~\ref{theorem90H}] 	Let $F$ be a field of characteristic different from $2$ such that $	(G_F(2),\theta_F)\in\mathcal C.$ Then Conjecture~\ref{claimh90} holds for $F$.
\end{theorem*}

The decomposition criterion also applies to the principal previously
known positive cases of the finite-dimensional H-conjecture.
Quasi-Pythagorean fields with finitely many square classes, for which
the conjecture was established by Kijima and Nishi \cite{iwakami, kijima1, kijima2}, are covered by
Ware's decomposition of their maximal pro-$2$ Galois groups
\cite[Corollary~1]{ware1}. Pre-$2$-henselian fields \cite{dario3} admit the required decomposition by Example~\ref{ex:pre-2-henselian-fields}, while the families constructed by Berman \cite{berman} and Kula \cite{kula1, kula2}, recalled in
Examples~\ref{ex:berman1}, \ref{ex:berman2}, and~\ref{ex:kula}, are
treated in Section~\ref{sec:revisiting}. Consequently,
Theorem~\ref{thm:H90-decomposition} provides a common
Galois-theoretic framework for these previously known cases.

Although the present paper concerns the Kaplansky radical, analogous
questions arise for odd primes. In \cite{dario1}, we studied
the Galois theory of fields admitting valuations compatible with the
$p$-radical, a natural generalization of the Kaplansky radical.
More recently, Blumer, Feuerpfeil, Lopes, and Quadrelli
\cite{blumer-quadrelli} introduced a cohomological analogue of the
Kaplansky radical for arbitrary profinite groups and primes $p$. For
absolute Galois groups, their construction recovers the classical
Kaplansky radical when $p=2$ and the $p$-radical for arbitrary $p$.
They also formulated a group-theoretic analogue of the H-conjecture and studied it for several classes of pro-$p$ groups.

The emphasis of the present paper is different. We work directly in
the quadratic setting, where the Kaplansky radical arises naturally,
and focus on free pro-$2$ product decompositions of maximal pro-$2$
Galois groups. Our main purpose is to show how such decompositions
determine the Kaplansky radical and yield a criterion for the
H-conjecture, which remains open in general. This approach also
clarifies the connections between the Kaplansky radical, quadratic
forms, rigid elements, and compatible valuations.



Throughout the paper, all fields are assumed to have characteristic
different from $2$. All subgroups of pro-$2$ groups are understood to
be closed, and all homomorphisms are continuous.


The paper is organized as follows.

\begin{enumerate}[\upshape$\bullet$]
	\item Section~\ref{sec:prelims} recalls the necessary background on
	the Kaplansky radical, quadratic forms, and maximal pro-$2$ Galois
	groups.
	
	\item In Section~\ref{sec:R(F)-and-fre-products}, we characterize
	fields whose maximal pro-$2$ Galois groups are free and determine the
	behavior of the Kaplansky radical under free pro-$2$ product
	decompositions.
	
	\item Section~\ref{sec:main-result} studies how such a decomposition
	behaves after passing to a quadratic extension. Using a pro-$2$ version of the Kurosh subgroup theorem, we prove
	Theorem~\ref{thm:H90-decomposition}.
	
	\item In Section~\ref{sec:the-elementary-case}, we recall the
	elementary type construction in the category of cyclotomic pro-$2$ pairs and 	determine the Kaplansky radical for fields realizing the basic pairs and cyclotomic semidirect products. Following the recursive construction of the elementary type class, we then prove
	Theorem~\ref{thm:claim2-ETG2} and obtain
	Theorem~\ref{theorem90H} as a consequence.
	
	\item Finally, Section~\ref{sec:revisiting} revisits the
	classical constructions of Berman and Kula and shows that they
	satisfy the decomposition criterion. We also explain how
	pre-$2$-henselian fields and quasi-Pythagorean fields with finitely
	many square classes fit into the same framework.
\end{enumerate}

\section{Preliminaries}\label{sec:prelims}

In this section, we collect the basic facts about the Kaplansky radical
needed in the sequel. Most of the material is classical and can be
found in \cite[Chapter XII, \S6]{lam1}. We begin with two equivalent
characterizations that will be used throughout the paper.

For a field $F$ and an element $a\in\dot F\setminus\dot F^2$, let
$K=F(\sqrt a)$, and let
\[
N=N_{K/F}:\dot K\longrightarrow\dot F
\]
be the norm map of the quadratic extension $K/F$. Define
\[
D_F\langle1,-a\rangle=N(\dot K),
\]
the image of the norm map $N$, or equivalently, the set of nonzero
values represented by the binary quadratic form $X^2-aY^2$. By \eqref{eq:radical-as-intersection}, the Kaplansky radical of $F$ is given by
\[
R(F)=\bigcap_{\alpha\in\dot F}D_F\langle1,\alpha\rangle.
\]
In particular, $R(F)$ is a subgroup of $\dot F$. Moreover, combining this definition with the equivalence 
\begin{equation}\label{eq:groups-norms}
a\in D_F\langle1,-b\rangle
\Longleftrightarrow
b\in D_F\langle1,-a\rangle,
\end{equation}
we obtain the equivalent description
\begin{equation}\label{eq:R(F)-characterization} R(F)=\Big\{c\in\dot F \ \ \bigl\vert \ \  D_F\langle1,-c\rangle=\dot F\Big\},\end{equation}
which is frequently taken as the definition of the Kaplansky radical in the literature.

A second characterization arises naturally in the context of central
simple algebras. Let $(F;a,b)$ denote the quaternion algebra over $F$
generated by elements $i$ and $j$ satisfying $i^2=a,$ $j^2=b$, $ij=-ji$, 
and let $\twotorsion$ denote the $2$-torsion subgroup of the Brauer
group of $F$. Since \eqref{eq:groups-norms} is also equivalent to $[(F;a,b)]=0$ in $\twotorsion$, we obtain 
\begin{equation}\label{eq:R(F)-characterization-quaternions}R(F)=\Big\{a\in\dot F  \ \ \bigl\vert \ \  [(F;a,b)]=0 \text{ for every }b\in\dot F\Big\}, 
\end{equation}
that is, $R(F)/\dot F^2$ is the radical of the symmetric bilinear
pairing
\[
\dot F/\dot F^2\times\dot F/\dot F^2
\longrightarrow
\twotorsion,
\]
defined by $(\bar a,\bar b)\longmapsto[(F;a,b)].$

From the definitions above, we obtain the following inclusions:
\begin{equation}\label{chain}
\dot F^2 \subseteq R(F)  \subseteq  D_F\langle1,1\rangle \subseteq \dot F.
\end{equation}
	
\begin{rem}\label{rem:chain}
The three possible equalities in (\ref{chain}) involving $R(F)$ lead to the following terminology.
\begin{enumerate}[\upshape(1)]
\item \textbf{(Trivial case)} 
If \[R(F)=\dot F^2,\] we say that $F$ has a \emph{trivial Kaplansky radical}. This is the case, for example, for Pythagorean fields, $2$-Demushkin fields (see Proposition~\ref{prop:demushkin-radical}), and the $2$-henselian fields
covered by Proposition~\ref{prop:radical-of-2hens}.
\item \textbf{(Quasi-Pythagorean case)} If 	\[R(F)=D_F\langle1,1\rangle\]		
then $F$ is called a \emph{quasi-Pythagorean field}. Such fields form
one of the principal classes of fields with nontrivial Kaplansky
radical and have been studied extensively \cite{ kijima1, kijima2, ware1}.
\item \textbf{(Free case)} If \[R(F)=\dot F,\] then \eqref{eq:R(F)-characterization} implies $D_F\langle1,-a\rangle=\dot F$, for every $a\in\dtF.$ Typical examples include finite fields and Laurent series fields $k((T))$, where $k$ is quadratically closed. The terminology is justified in
Section~\ref{sec:R(F)-and-fre-products}, where we prove that this
condition is equivalent to $G_F(2)$ being a free pro-$2$ group.
\end{enumerate}
\end{rem}

\section{Kaplansky radical and free pro-$2$ products}
\label{sec:R(F)-and-fre-products}

The decomposition criterion of
Theorem~\ref{thm:H90-decomposition} relies on the behavior of the
Kaplansky radical under free pro-$2$ products. In this section, we
first characterize the fields whose maximal pro-$2$ Galois groups are
free and then determine the Kaplansky radical from such a
decomposition.
	

We begin with the free case. Recall that a free pro-$2$ group is characterized by its universal property; see \cite[Definition~3.5.14]{neukirch}.

\begin{prop}\label{prop:free-radical} For a field $F$, the following conditions are equivalent. 
\begin{enumerate}[\upshape(i)] 
			\item $R(F)=\dot F$; 
			\item $\twotorsion=\{0\}$; 
			\item $G_F(2)$ is a free pro-$2$ group. 
\end{enumerate} \end{prop} 
\begin{proof} By \eqref{eq:R(F)-characterization-quaternions}, $R(F)=\dot F$ if and only if $[(F;a,b)]=0$ for every $a,b\in\dot F$. Since quaternion algebras generate the $2$-torsion part of the Brauer group by the Merkurjev theorem, this is equivalent to $\twotorsion=\{0\}$. The equivalence between $\twotorsion=\{0\}$ and the freeness of $G_F(2)$ follows from \cite[Corollary~3.8 and Corollary~3.2]{ribes1}. 
\end{proof} 
	
Proposition~\ref{prop:free-radical} is the $2$-analogue of a more general result appearing in \cite{dario3}.

\begin{ex}\label{ex:rank1} 	The free pro-$2$ group of rank one is isomorphic to
$\mathbb Z_2$, 	the additive group of the $2$-adic integers. It occurs as the maximal pro-$2$ Galois group of finite fields and of $\mathbb C((T))$ 	\cite[p.~392]{ware3}.
\end{ex}


We now turn to free pro-$2$ products, the coproducts in the category of
pro-$2$ groups; see \cite[Section 4.1]{ribes2}. To formulate
the criterion used below, we briefly recall the relevant cohomological
notation.

If $G$ is a profinite group and $M$ is a discrete $G$-module, we denote by $H^n(G,M)$ its $n$-th continuous cohomology group, defined using continuous cochains; see \cite[Chapter~I, \S2]{neukirch}. Throughout
this paper, the coefficient module is $\mathbb F_2$ with trivial action. For the maximal pro-$2$ Galois group of a field $F$, Kummer theory and the
norm-residue isomorphism give, respectively,
\begin{equation}\label{eq:kummer-norm-residue}
H^1(G_F(2),\mathbb F_2) \cong \dot F/\dot F^2,	\qquad 	H^2(G_F(2),\mathbb F_2)	\cong
\twotorsion.
\end{equation}

Under these identifications, the cup product
\[\cup:H^1(G_F(2),\mathbb F_2)\times H^1(G_F(2),\mathbb F_2)\longrightarrow H^2(G_F(2),\mathbb F_2) \]
corresponds to the quaternion pairing $(\bar a,\bar b)\longmapsto[(F;a,b)].$

The following criterion therefore translates a free pro-$2$ product
decomposition of $G_F(2)$ into conditions on square classes and
quaternion algebras.

\begin{thm}\label{teo:decomposition-neu-traduzido} Let $F$ be a field and let $F_1,\ldots,F_n$ be intermediate fields of $F(2)/F$. Then \[ G_F(2) \cong G_{F_1}(2)\ast_2\cdots\ast_2G_{F_n}(2) \] if and only if the following conditions hold: 
\begin{enumerate}[\upshape(i)]
			\item The inclusions $F\hookrightarrow F_i$ induce an isomorphism \[ \dot F/\dot F^2 \longrightarrow \prod_{i=1}^{n}\dot F_i/\dot F_i^2. \] 
			\item The scalar extension maps induce a monomorphism \[ \twotorsion \longrightarrow \prod_{i=1}^{n}{}_2\Br{F_i}. \] 
\end{enumerate} \end{thm} 
\begin{proof} By Neukirch's cohomological criterion \cite[Satze~4.2 and~4.3]{neukirch2}, together with the Kummer and norm-residue isomorphisms	\eqref{eq:kummer-norm-residue}, the asserted decomposition is equivalent	to the restriction maps
	\[ 	H^i(G_F(2),\mathbb F_2)\longrightarrow	\bigoplus_{j=1}^{n}H^i(G_{F_j}(2),\mathbb F_2),	\qquad i=1,2,	\]
		being an isomorphism for $i=1$ and a monomorphism for $i=2$. Under (\ref{eq:kummer-norm-residue}), these correspond precisely to
		conditions~\textup{(i)} and~\textup{(ii)}, respectively.
\end{proof} 
	
The preceding criterion gives the following description of the Kaplansky radical under free pro-$2$ products. 
	
\begin{thm}\label{radical-of-fre-product} Let $F$ be a field and let $F_1,\ldots,F_n$ be intermediate fields of $F(2)/F$ such that \[ G_F(2) \cong G_{F_1}(2)\ast_2\cdots\ast_2G_{F_n}(2). \] 
Then the following statements hold:
\begin{enumerate}[\upshape(a)]
	\item For every $x\in\dot F$, \[ D_F\langle1,x\rangle = \dot F\cap \bigcap_{i=1}^{n}D_{F_i}\langle1,x\rangle. \] 
	\item The Kaplansky radical of $F$ is given by \[ R(F) = \dot F\cap \bigcap_{i=1}^{n}R(F_i). \] 
	\item If $R(F_i)=\dot F_i^{\,2}$ for every $i$, then $R(F)=\dot F^{\,2}.$ 	Likewise, if $R(F_i)=\dot F_i$ for every $i$, then $R(F)=\dot F.$
	\end{enumerate} 
\end{thm} 
\begin{proof}  
\noindent\textup{(a)} Scalar extension gives the inclusion from left to
right. Conversely, let \[ y\in \dot F\cap \bigcap_{i=1}^{n}D_{F_i}\langle1,x\rangle. \] By \eqref{eq:groups-norms} and \eqref{eq:R(F)-characterization-quaternions}, the image of the class $[(F;-x,y)]\in \twotorsion$ is trivial in ${}_2\Br{F_i}$ for every $i$. 
Since the map in Theorem~\ref{teo:decomposition-neu-traduzido}\textup{(ii)} is injective, it follows that $[(F;-x,y)]=0$. Therefore, $y\in D_F\langle1,x\rangle$. 
\\
\noindent\textup{(b)}
By Theorem~\ref{teo:decomposition-neu-traduzido}\textup{(i)}, every
square class of $F_i$ is represented by an element of $\dot F$. Since
$D_{F_i}\langle1,x\rangle$ depends only on the square class of $x$, we
have
\[
\bigcap_{x\in\dot F}D_{F_i}\langle1,x\rangle=R(F_i)
\]
for every $i$. The result now follows by taking intersections over all
$x\in\dot F$ in the equality of part~\textup{(a)}.
\\
\noindent\textup{(c)}
Suppose first that $R(F_i)=\dot F_i^2$ for every $i$. By part~\textup{(b)},
\[R(F)=\dot F\cap\bigcap_{i=1}^{n}\dot F_i^2.\]
The injectivity of the map in Theorem~\ref{teo:decomposition-neu-traduzido}\textup{(i)} then implies
that $R(F)=\dot F^2.$ If instead $R(F_i)=\dot F_i$ for every $i$, then part~\textup{(b)} immediately gives $R(F)=\dot F.$
\end{proof}

The following special case, consisting of one free factor and factors
with trivial Kaplansky radical, will be used repeatedly.
	
\begin{cor}\label{cor:free-factor-F_1} Let $F$ be a field and let $F_1,\ldots,F_n$ be intermediate fields of $F(2)/F$ such that \[	G_F(2)	\cong	G_{F_1}(2)\ast_2\cdots\ast_2G_{F_n}(2).	\]
Assume that $G_{F_1}(2)$ is free pro-$2$ and that $R(F_i)=\dot F_i^2$	for every $i\geq2$. Then:
	\begin{enumerate}[\upshape(a)]
		\item The Kaplansky radical of $F$ is given by
		\[	R(F)		=		\dot F\cap\bigcap_{i=2}^{n}\dot F_i^2.		\]
	\item The inclusion $F\hookrightarrow F_1$ induces an isomorphism 
	\[	R(F)/\dot F^2	\longrightarrow	\dot F_1/\dot F_1^2.\]
	Consequently,
\[	\dim_{\mathbb F_2}\bigl(R(F)/\dot F^2\bigr)	=	\dim_{\mathbb F_2}\bigl(\dot F_1/\dot F_1^2\bigr)=	\operatorname{rk}\bigl(G_{F_1}(2)\bigr).\]
\end{enumerate} 
\end{cor} 
\begin{proof} \noindent\textup{(a)} By Proposition~\ref{prop:free-radical}, $R(F_1)=\dot F_1.$ 
The result therefore follows immediately from Theorem~\ref{radical-of-fre-product}\textup{(b)}.
\\
\noindent\textup{(b)}
Let
\[\Phi:\dot F/\dot F^2\longrightarrow\prod_{i=1}^{n}\dot F_i/\dot F_i^2\]
be the isomorphism of Theorem~\ref{teo:decomposition-neu-traduzido}\textup{(i)}. 
By part~\textup{(a)}, a class in $\dot F/\dot F^2$ belongs to
$R(F)/\dot F^2$ if and only if all components of its image under
$\Phi$, except possibly the first, are trivial. Since $\Phi$ is an
isomorphism, projection onto the first factor restricts to the required isomorphism.  Finally, the equality of dimensions follows immediately, while the last equality is a consequence of Kummer theory; see (\ref{eq:kummer-norm-residue}).
\end{proof}
	
\section{The Main Result}
\label{sec:main-result}

In this section, we prove Theorem~\ref{thm:H90-decomposition}. Let $F$ be a field and let $F_0$ and $F_1$ be intermediate fields of $F(2)/F$ such that \begin{equation}\label{eq:decomposition}
G_F(2) \cong G_{F_0}(2)\ast_2G_{F_1}(2),
\end{equation}
where $G_{F_0}(2)$ is free pro-$2$ and $F_1$ has trivial Kaplansky
radical. Let
\[K=F(\sqrt a), \qquad a\in\dot F\setminus\dot F^2.\]

Our first step is to determine how the decomposition
\eqref{eq:decomposition} behaves after passing to the index-two
subgroup $G_K(2)$ and to use the resulting decomposition to describe
$R(K)$. The main tool is the pro-$2$ Kurosh subgroup theorem
(Theorem~\ref{thm:kurosh}), which is a special case of
\cite[Theorem~4.2.1, p.~208]{neukirch}; see also \cite{bnw}.

\begin{thm}[Pro-$2$ Kurosh subgroup theorem] 	\label{thm:kurosh}
Let $G$ be a pro-$2$ group with 
\[	G\cong G_0\ast_2G_1,	\]
where $G_0$ and $G_1$ are closed subgroups of $G$, and let $H$ be an 	open subgroup of $G$. Choose double-coset decompositions
\[	G=\coprod_{i=1}^{n}G_0a_iH	\qquad\text{and}\qquad	G=\coprod_{j=1}^{m}G_1b_jH,	\]
where $a_i, b_j \in G$ for all $i, j$. Then
\[	H\cong	\left(	\mathop{\ast_2}_{i=1}^{n}	H\cap a_i^{-1}G_0a_i	\right)
\ast_2	\left(	\mathop{\ast_2}_{j=1}^{m}	H\cap b_j^{-1}G_1b_j	\right)	\ast_2\mathcal F,	\]
where $\mathcal F$ is a free pro-$2$ group of rank
\[	\operatorname{rk}(\mathcal F)	=	1+(G:H)-n-m.	\]
\end{thm}
We apply Theorem~\ref{thm:kurosh} to the decomposition
\eqref{eq:decomposition} and to the index-two subgroup
$H = G_K(2)\subseteq G_F(2).$ By Galois correspondence,
\[G_{F_i}(2)\subseteq G_K(2)\quad\Longleftrightarrow\quad K\subseteq F_i\quad\Longleftrightarrow\quad a\in\dot F_i^2,\]
for $i=1,2$. Moreover, Corollary~\ref{cor:free-factor-F_1}\textup{(a)} gives $R(F)=\dot F\cap\dot F_1^2,$ and hence
\begin{equation}\label{eq:a in R(F) or not}
a\in R(F) \quad\Longleftrightarrow\quad a\in\dot F_1^2.
\end{equation}

Thus the double cosets contributed by $G_{F_1}(2)$ are determined by
whether $a$ belongs to $R(F)$. The following theorem describes the
corresponding decompositions of $G_K(2)$ and the resulting
descriptions of $R(K)$.

\begin{thm}\label{thm:decomposition}
	Let $F$ be a field of characteristic different from $2$ satisfying
	\eqref{eq:decomposition}, and let 
	\[	K=F(\sqrt a),\qquad a\in\dot F\setminus\dot F^2.	\]
 Then the following statements hold.
	\begin{enumerate}[\upshape(i)]
		\item Suppose that $a\notin R(F)$, and set $K_1=F_1(\sqrt a).$
		Then there exists an intermediate field $L_0$ of $F(2)/K$ such that
		$G_{L_0}(2)$ is free pro-$2$ and
		\[
		G_K(2)\cong G_{L_0}(2)\ast_2G_{K_1}(2).
		\]
		Consequently, $		R(K)=\dot K\cap R(K_1).$		
		\item Suppose that $a\in R(F)$. Choose
		$\sigma\in G_F(2)\setminus G_K(2)$ and set $F_1^\sigma=\sigma(F_1).$ Then there exists an intermediate field $L_0$ of $F(2)/K$ such that
		$G_{L_0}(2)$ is free pro-$2$ and
		\[
		G_K(2)\cong
		G_{L_0}(2)\ast_2G_{F_1}(2)\ast_2G_{F_1^\sigma}(2).
		\]
		Moreover, $R(K)=\dot K\cap\dot F_1^2\cap\dot{F_1^\sigma}^2.$
	\end{enumerate}
\end{thm}
\begin{proof}
\noindent\textup{(i)} Assume that $a\notin R(F)$. By \eqref{eq:a in R(F) or not}, we have $a\notin\dot F_1^2.$ Set $K_i=F_i(\sqrt a)$ for $i=0,1$. In particular,
\[
G_{F_1}(2)\nsubseteq G_K(2)
\qquad\text{and}\qquad
G_{F_1}(2)\cap G_K(2)=G_{K_1}(2).
\]
We distinguish two cases according to whether $a$ is a square in
$F_0$.

\medskip

\noindent\textit{Case 1: $a\notin\dot F_0^2$.} In this case, the relevant fields are arranged as follows:
\[
\xymatrix@C=4.0em@R=2.3em{
	& F(2) & \\
	K_0=F_0(\sqrt a)\ar@{-}[ur]
	&&
	K_1=F_1(\sqrt a)\ar@{-}[ul]
	\\
	F_0\ar@{-}[u]
	&
	K=F(\sqrt a)\ar@{-}[uu]\ar@{-}[ul]\ar@{-}[ur]
	&
	F_1\ar@{-}[u]
	\\
	&
	F\ar@{-}[u]\ar@{-}[ul]\ar@{-}[ur]
	&
}\]
Since $a\notin\dot F_0^2$ and $a\notin\dot F_1^2,$ neither $G_{F_0}(2)$ nor $G_{F_1}(2)$ is contained in $G_K(2)$. Moreover,
\[
G_{F_i}(2)\cap G_K(2)=G_{K_i}(2),
\qquad i=0,1.
\]
Since $G_K(2)$ is normal of index two in $G_F(2)$, each factor
contributes one double coset. Therefore,
Theorem~\ref{thm:kurosh} gives
\[
G_K(2)\cong
G_{K_0}(2)\ast_2G_{K_1}(2)\ast_2\mathcal F,
\]
where $\operatorname{rk}(\mathcal F)=1+2-1-1=1.$ Thus $\mathcal F\cong\mathbb Z_2.$  Since $G_{K_0}(2)$ is a closed subgroup of the free pro-$2$ group $G_{F_0}(2)$, it is free pro-$2$. Consequently, the group
$\mathcal G_0 = G_{K_0}(2)\ast_2\mathcal F$ 
is free pro-$2$, and
\[G_K(2)\cong\mathcal G_0\ast_2G_{K_1}(2).\]
\noindent\textit{Case 2: $a\in\dot F_0^2$.} In this case, $K_0=F_0$ and the corresponding fields are arranged as follows:
\[
\xymatrix@C=3.2em@R=2.3em{
	& F(2) & \\
	K_0 = F_0\ar@{-}[ur]
	&
	F_0^\tau\ar@{-}[u]
	&
	K_1=F_1(\sqrt a)\ar@{-}[ul]
	\\
	&
	K=F(\sqrt a)\ar@{-}[u]\ar@{-}[ul]\ar@{-}[ur]
	&
	F_1\ar@{-}[u]
	\\
	&
	F\ar@{-}[u]\ar@{-}[ur]
	&
}
\]
We have $G_{F_0}(2)\subseteq G_K(2),$ $G_{F_1}(2)\nsubseteq G_K(2),$ and
$G_{F_1}(2)\cap G_K(2)=G_{K_1}(2).$ Choose $\tau\in G_F(2)\setminus G_K(2).$ Since $G_K(2)$ is normal of index two in $G_F(2)$, the factor
$G_{F_0}(2)$ contributes two double cosets, represented by $1$ and
$\tau$, whereas $G_{F_1}(2)$ contributes one. The free factor in the
Kurosh decomposition has rank $1+2-2-1=0.$ Consequently,
\[G_K(2)\cong G_{F_0}(2) \ast_2 \tau^{-1}G_{F_0}(2)\tau \ast_2 G_{K_1}(2).\]
Both $G_{F_0}(2)$ and $\tau^{-1}G_{F_0}(2)\tau$ are free pro-$2$ groups. Therefore, the group $\mathcal G_0 = G_{F_0}(2) \ast_2 \tau^{-1}G_{F_0}(2)\tau $ is free pro-$2$, and
\[G_K(2)\cong\mathcal G_0\ast_2G_{K_1}(2).\]
Thus, in either case, there exists a free pro-$2$ subgroup
$\mathcal G_0$ of $G_K(2)$ such that $G_K(2)\cong\mathcal G_0\ast_2G_{K_1}(2).$ Let $L_0$ be the fixed field of $\mathcal G_0$ in $F(2)$. Then $L_0$ is intermediate in $F(2)/K$, $G_{L_0}(2)=\mathcal G_0,$ and therefore
\[
G_K(2)\cong G_{L_0}(2)\ast_2G_{K_1}(2),
\]
where $G_{L_0}(2)$ is free pro-$2$. By Proposition~\ref{prop:free-radical}, $R(L_0)=\dot L_0.$ Finally, it follows from Theorem~\ref{radical-of-fre-product}\textup{(b)} that
\[
R(K)=\dot K\cap R(K_1).
\]
\noindent\textup{(ii)}
Suppose that $a\in R(F)$. By
\eqref{eq:a in R(F) or not}, we have $a\in\dot F_1^2.$  
Consequently, $K\subseteq F_1$ and $G_{F_1}(2)\subseteq G_K(2).$ 
On the other hand, $a\notin\dot F_0^2.$  Indeed, if $a$ were a square in both $F_0$ and $F_1$, the injectivity in Theorem~\ref{teo:decomposition-neu-traduzido}\textup{(i)} would
imply that $a\in\dot F^2,$ contrary to the choice of $a$. Set $K_0=F_0(\sqrt a).$ Choose $\sigma\in G_F(2)\setminus G_K(2)$ and set 
$F_1^\sigma=\sigma(F_1).$ Since $\sigma(K)=K$, we also have
$K\subseteq F_1^\sigma.$ In particular, $F_1(\sqrt a)=F_1$ and $F_1^\sigma(\sqrt a)=F_1^\sigma.$ The relevant fields are arranged as follows:
\[
\xymatrix@C=3.2em@R=2.3em{
	& F(2) & \\
	K_0=F_0(\sqrt a)\ar@{-}[ur]
	&
	F_1^\sigma\ar@{-}[u]
	&
	F_1\ar@{-}[ul]
	\\
	F_0\ar@{-}[u]
	&
	K=F(\sqrt a)\ar@{-}[u]\ar@{-}[ul]\ar@{-}[ur]
	&
	\\
	F\ar@{-}[u]\ar@{-}[ur]
	&&
}\]
Since $G_K(2)$ is normal of index two in $G_F(2)$, the factor
$G_{F_0}(2)$ contributes one double coset, whereas
$G_{F_1}(2)$ contributes two, represented by $1$ and
$\sigma^{-1}$. The free factor in the Kurosh decomposition has rank
$1+2-1-2=0.$ Therefore,
\[G_K(2)\cong G_{K_0}(2)  \ast_2 G_{F_1}(2) \ast_2 G_{F_1^\sigma}(2). \]
Since $G_{K_0}(2)$ is a closed subgroup of the free pro-$2$ group
$G_{F_0}(2)$, it is free pro-$2$. Taking $L_0=K_0=F_0(\sqrt a),$
we obtain 
\[G_K(2)\cong G_{L_0}(2)  \ast_2 G_{F_1}(2) \ast_2 G_{F_1^\sigma}(2). \]
Finally, $\sigma$ induces an isomorphism $F_1\cong F_1^\sigma.$ Since $F_1$ has trivial Kaplansky radical, so does $F_1^\sigma$. 
Corollary~\ref{cor:free-factor-F_1}\textup{(a)} now gives
\[
R(K)
=
\dot K
\cap
\dot F_1^2
\cap
\dot{F_1^\sigma}^2.
\]
\end{proof}

The decomposition obtained in Theorem~\ref{thm:decomposition} provides
the key ingredient for the proof of the main result. 

\begin{thm}\label{thm:H90-decomposition}
	Let $F$ be a field of characteristic different from $2$, and let
	$F_0,F_1$ be intermediate fields of $F(2)/F$ such that
	\[
	G_F(2)\cong G_{F_0}(2)\ast_2G_{F_1}(2),
	\]
	where $G_{F_0}(2)$ is free pro-$2$ and $F_1$ has trivial Kaplansky
	radical. Let $K=F(\sqrt a)$, where
	$a\in\dot F\setminus\dot F^2$, and let $N=N_{K/F}:\dot K\longrightarrow\dot F$
	be the norm map. Then
	\[
	N^{-1}(R(F))=\dot FR(K).
	\]
\end{thm}
\begin{proof} By the Knebusch norm principle, $	N(R(K))\subseteq R(F)$.
Hence, for every $c\in\dot F$ and $r\in R(K)$, $N(cr)=c^2N(r)\in R(F),$
because $\dot F^2\subseteq R(F)$. Therefore, $\dot F\,R(K)\subseteq N^{-1}(R(F)).$ 
For the reverse inclusion, let $x\in\dot K$ satisfy $N(x)\in R(F)$. 
By Corollary~\ref{cor:free-factor-F_1}\textup{(a)},
\[R(F)=\dot F\cap\dot F_1^2.\]
In particular, $N(x)\in\dot F_1^2$. We distinguish the two cases appearing
in Theorem~\ref{thm:decomposition}.
\medskip
\\ 
\noindent\textit{Case 1: $a\notin R(F)$.}
Set $K_1=F_1(\sqrt a)$. By \eqref{eq:a in R(F) or not},
$a\notin\dot F_1^2,$ so $K_1/F_1$ is a quadratic extension. Let
\[
N_1=N_{K_1/F_1}:\dot K_1\longrightarrow\dot F_1
\]
be its norm map. Since $K\subseteq K_1$ and the nontrivial automorphism of $K_1/F_1$ restricts to the nontrivial automorphism of $K/F$, we have
\[N_1(x)=N(x)\in\dot F_1^2.\] It follows from Hilbert's Theorem~90, in the form \eqref{eq:H90-quadratic}, that $x\in\dot F_1\dot K_1^2.$ By Theorem~\ref{teo:decomposition-neu-traduzido}\textup{(i)}, $\dot F_1=\dot F\dot F_1^2.$  Hence
\[x\in\dot F\dot K_1^2\subseteq\dot F R(K_1).\]
Thus $x=\alpha r$ for some $\alpha \in\dot F$ and $r\in R(K_1)$. Since
$x,\alpha \in\dot K$, we also have $r\in\dot K$. Therefore,
Theorem~\ref{thm:decomposition}\textup{(i)} gives
\[
r\in\dot K\cap R(K_1)=R(K),
\]
and consequently $x\in\dot F R(K)$.
\medskip \\
\noindent\textit{Case 2: $a\in R(F)$.} By \eqref{eq:a in R(F) or not},
$a\in\dot F_1^2,$ and hence $K\subseteq F_1.$ Again, Theorem~\ref{teo:decomposition-neu-traduzido}\textup{(i)}
gives $\dot F_1=\dot F\dot F_1^2.$ Thus $x=\alpha t^2$ 
for some $\alpha\in\dot F$ and $t\in\dot F_1$. Choose
$\sigma\in G_F(2)\setminus G_K(2)$ as in Theorem~\ref{thm:decomposition}\textup{(ii)}, and set $F_1^\sigma=\sigma(F_1).$ Since $N(x)\in R(F)\subseteq\dot F_1^2$ 
and $N(x)\in\dot F$ is fixed by $\sigma$, we also have $N(x)\in\dot{F_1^\sigma}^2.$ Moreover,
\[N(x)=x\sigma(x)=\alpha^2t^2\sigma(t)^2.\]
Both $\alpha^2$ and $\sigma(t)^2$ 
belong to $\dot{F_1^\sigma}^2$. it follows that
$t^2\in(\dot F_1^\sigma)^2$. On the other hand, $t^2=\alpha^{-1}x\in\dot K.$ Therefore,
\[t^2\in\dot K\cap\dot F_1^2\cap\dot{F_1^\sigma}^2 = R(K) \]
by Theorem~\ref{thm:decomposition}\textup{(ii)}. Hence $x=\alpha t^2\in\dot F\,R(K),$ which completes the proof.
\end{proof}

Whether the converse of Theorem~\ref{thm:H90-decomposition} holds remains open. Nevertheless, we shall see that the families considered in this paper satisfy both its decomposition hypothesis and its norm identity,
providing further evidence that these two properties are closely
related.

\section{The Elementary Type Case}
\label{sec:the-elementary-case}

The decomposition criterion established in the preceding section was
formulated entirely in terms of maximal pro-$2$ Galois groups, since
the decompositions considered there were naturally compatible with
their canonical cyclotomic characters and there was no need to specify
them explicitly. To apply this criterion to the elementary type
construction, however, we must work in the category of cyclotomic
pro-$2$ pairs, since the cyclotomic character determines the
semidirect product operation used in that construction.

Recall that the elementary type class $\mathcal C$ is the class of
realizable cyclotomic pro-$2$ pairs generated from the basic pairs
associated with finitely generated free pro-$2$ groups and
$2$-Demushkin fields by means of free pro-$2$ products and cyclotomic
semidirect products. The pro-$2$ version of the Elementary Type
Conjecture predicts that every realizable cyclotomic pair whose
underlying pro-$2$ group is finitely generated belongs
to $\mathcal C$.
	
The purpose of this section is to show that the decomposition required
by Theorem~\ref{thm:H90-decomposition} is compatible with the
recursive construction of the elementary type class $\mathcal C$.
More precisely, we prove that, whenever
\[
(G_F(2),\theta_F)\in\mathcal C,
\]
there exist intermediate fields $L,H$ of $F(2)/F$ such that, in the
category of cyclotomic pro-$2$ pairs,
\[
(G_F(2),\theta_F)
\cong
(G_L(2),\theta_L)\ast_2(G_H(2),\theta_H),
\]
where $G_L(2)$ is free pro-$2$ and $R(H)=\dot H^2.$ The corresponding decomposition of the underlying pro-$2$ groups satisfies the hypotheses of
Theorem~\ref{thm:H90-decomposition}. Consequently, the H-conjecture
holds in the elementary type case.




The section is organized as follows.

\begin{enumerate}[\upshape$\bullet$]
\item In Subsection~\ref{subsec:R(F)-and-Demushkin}, we determine the
Kaplansky radical for fields realizing the two basic families of the
elementary type construction. The free case follows from
Proposition~\ref{prop:free-radical}, while the $2$-Demushkin case
follows from the nondegeneracy of the cup product.
	
\item In Subsection~\ref{subsec:ETG}, we recall the elementary type
construction in the category of cyclotomic pro-$2$ pairs, including
its basic pairs and the two operations used in its recursive
definition.


\item In Subsection~\ref{sec:semi-direct}, we study fields whose
canonical cyclotomic pro-$2$ pairs arise from cyclotomic semidirect products. Using results on rigid elements and $2$-henselian valuations, we prove that these fields have trivial Kaplansky radical.

\item Finally, in Subsection~\ref{sec:radical-ETG}, we combine these results with the behavior of the Kaplansky radical under free
	pro-$2$ products. Following the recursive construction of
	$\mathcal C$,  we prove Theorem~\ref{thm:claim2-ETG2} and derive
	Theorem~\ref{theorem90H}.
\end{enumerate}

\subsection{The Kaplansky radical and $2$-Demushkin Galois groups}
\label{subsec:R(F)-and-Demushkin}
Finitely generated free pro-$2$ groups and $2$-Demushkin groups underlie the two basic families of cyclotomic pro-$2$ pairs used in the elementary type construction. The free case was treated in Proposition~\ref{prop:free-radical}. We now recall the corresponding result for $2$-Demushkin fields, which follows directly from the nondegeneracy of the cup product.

\begin{defn} \cite{demushkin1, demushkin2} \label{def:demushkin-group}
	Let $p$ be a prime. A pro-$p$ group $G$ is called a
	\emph{Demushkin group}if:
	\begin{enumerate}[\upshape(a)]
		\item $H^1(G,\mathbb F_p)$ is finite-dimensional;
		\item $H^2(G,\mathbb F_p)\cong\mathbb F_p$;
		\item the cup product, denoted by
		\[
		\cup:
		H^1(G,\mathbb F_p)\times H^1(G,\mathbb F_p)
		\longrightarrow H^2(G,\mathbb F_p),
		\]
		is a nondegenerate bilinear pairing.
	\end{enumerate}
\end{defn}

Apart from the cyclic group of order $2$, every Demushkin group is
infinite. The classification of Demushkin groups  is due to Serre and Labute; see \cite[Theorems~3.9.11 and~3.9.19]{neukirch}. In what follows, we restrict to the case $p=2$.

\begin{defn}\label{def:demushkin-field-and-group}
	A field $F$ is called a \emph{$2$-Demushkin field} if $G_F(2)$ is a
	Demushkin group. A realizable Demushkin pro-$2$ group, that is, a
	Demushkin group isomorphic to $G_F(2)$ for some field $F$, will be
	called a \emph{$2$-Demushkin Galois group}.
\end{defn}

\begin{ex}\label{ex-demushkin-fields}
If $L/\Q_2$ is 	a finite extension, then $L$ is a $2$-Demushkin field, and the minimal number of generators of $G_L(2)$ is $[L:\Q_2]+2$. For $\Q_2$ itself, 
\[G_{\mathbb Q_2}(2)\cong\left\langle \alpha,\beta,\gamma \mathrel{\Big|} 
\alpha^2\beta^4[\beta,\gamma]=1 \right\rangle, \]
 where 
$[\beta,\gamma]$ denotes the commutator; see \cite[Table~5.2]{ware3}.
\end{ex}
	
Under the Kummer and norm-residue isomorphisms \eqref{eq:kummer-norm-residue}, the cup product corresponds to the quaternion pairing
\[
(\bar a,\bar b)\longmapsto[(F;a,b)]
\]
introduced in Section~\ref{sec:prelims}. Consequently, the radical of the cup-product pairing is identified with $R(F)/\dot F^2$. 

\begin{prop}\label{prop:demushkin-radical} 	If $F$ is a $2$-Demushkin field, then
$R(F)=\dot F^2.$
\end{prop}
\begin{proof} By \eqref{eq:R(F)-characterization-quaternions} and
\eqref{eq:kummer-norm-residue}, an element $a\in\dot F$ belongs to
$R(F)$ if and only if $(a)\cup(b)=0$ for every $b\in\dot F$, where $(a)$ and $(b)$ denote the corresponding classes in $H^1(G_F(2),\mathbb F_2)$. By the nondegeneracy of the cup product, the class $(a)$ is trivial. Hence $a\in\dot F^2$, and therefore $R(F)=\dot F^2$.
\end{proof}


\subsection{Elementary type cyclotomic pairs}
\label{subsec:ETG}

We now recall the elementary type construction in the category of
cyclotomic pro-$2$ pairs. The definitions below follow
\cite[Section~4]{ware3}; see also \cite{efrat-livro}. We include only the
material needed for the proof of Theorem~\ref{thm:claim2-ETG2}.


\begin{defn}\label{def:par-ciclotomico} 
A \emph{cyclotomic pro-$2$ pair} is a pair $(G,\theta)$, where $G$ is a 	pro-$2$ group and \[\theta: G \longrightarrow \mathbb Z_2^{\times}\]
is a continuous homomorphism. 
A {\it morphism} \[
\varphi:(G_1,\theta_1)\longrightarrow(G_2,\theta_2)
\]
of cyclotomic pairs is a continuous homomorphism
$\varphi:G_1\to G_2$ such that $\theta_1=\theta_2\circ\varphi.$ 
An {\it isomorphism} of cyclotomic pairs is a morphism whose underlying group homomorphism is an isomorphism.
Finally, a cyclotomic pair is called 	\emph{realizable} if it is isomorphic to \[(G_F(2),\theta_F)\] for some 
field $F$,  where $\theta_F$ is the canonical cyclotomic
character describing the action of $G_F(2)$ on the $2$-primary roots
of unity $\mu_{2^\infty}\subseteq F(2).$
\end{defn}

The elementary type construction starts from the following two families
of basic realizable cyclotomic pro-$2$ pairs. To streamline the
presentation, we describe these families directly in the realizable
setting, without reproducing the more detailed classification of the
underlying groups given in \cite{ware3}.

\begin{enumerate}[\upshape(i)]
	\item canonical cyclotomic pro-$2$ pairs
	\[
	(G_F(2),\theta_F),
	\]
	where $G_F(2)$ is a finitely generated free pro-$2$ group. This
	includes the rank-$0$ case, corresponding to the trivial pair $
	(\{1\},1);$
	
	\item canonical cyclotomic pro-$2$ pairs
	\[
	(G_F(2),\theta_F),
	\]
	where $F$ is a $2$-Demushkin field.  This family includes the
	exceptional finite case $G_F(2)\cong\mathbb Z/2.$
	More precisely, for a real closed field $F$, the nontrivial element
	$\sigma\in G_F(2)$ satisfies $\theta_F(\sigma)=-1.$
\end{enumerate}

In both families, the character is the canonical cyclotomic character
associated with the realizing field; see \cite{ware3}. Starting from these basic pairs, the elementary type construction repeatedly applies two operations: free pro-$2$ products and cyclotomic semidirect products. We now describe these operations, beginning with free pro-$2$ products.
\\

\paragraph{{\bf Free pro-$2$ product.}}
Let $(G_1,\theta_1)$ and $(G_2,\theta_2)$ be cyclotomic pro-$2$ pairs. Their 
free pro-$2$ product is defined by
\[
(G_1,\theta_1)\ast_2(G_2,\theta_2)
=
(G_1\ast_2G_2,\theta),
\]
where $\theta:G_1\ast_2G_2\longrightarrow \Z_2^{\times}$ is the unique continuous homomorphism satisfying 
\[
\left.\theta\right|_{G_i}=\theta_i,
\qquad i=1,2.
\]
Existence and uniqueness of $\theta$ follow from the universal property of the free pro-$2$ product.  
\\

\paragraph{\textbf{Cyclotomic semidirect product.}}
Let $(G,\theta)$ be a cyclotomic pro-$2$ pair and let $m\geq1$. The
group $G$ acts on $\mathbb Z_2^m$ by scalar multiplication through
$\theta$:
\[
g\cdot z=\theta(g)z,
\qquad
g\in G,\quad z\in\mathbb Z_2^m.
\]
The associated cyclotomic semidirect product is the pair
\[
\bigl(\mathbb Z_2^m\rtimes_\theta G,\widetilde\theta\bigr),
\]
where
\[
\widetilde\theta:
\mathbb Z_2^m\rtimes_\theta G
\longrightarrow
\mathbb Z_2^\times,
\qquad
\widetilde\theta(z,g)=\theta(g).
\]
Thus, $\widetilde\theta$ is trivial on the normal subgroup
$\mathbb Z_2^m$ and restricts to $\theta$ on $G$.



The realizability results recalled in \cite[Section~4]{ware3} show that
both free pro-$2$ products and cyclotomic semidirect products preserve
realizability. We may therefore define the elementary type class
within the category of realizable cyclotomic pro-$2$ pairs.

\begin{defn}\label{def:elementary-type-class}
	Let $\mathcal C$ denote the smallest class of realizable cyclotomic pro-$2$ pairs satisfying the following conditions:
	
	\begin{enumerate}[\upshape(i)]
		\item If $(G_F(2),\theta_F)$ is a realizable cyclotomic pair and
		$G_F(2)$ is a finitely generated free pro-$2$ group, then
		\[
		(G_F(2),\theta_F)\in\mathcal C.
		\]
		
		\item If $F$ is a $2$-Demushkin field, then
		\[
		(G_F(2),\theta_F)\in\mathcal C.
		\]
		
		\item If $(G_1,\theta_1),(G_2,\theta_2)\in\mathcal C$, then
		\[
		(G_1,\theta_1)\ast_2(G_2,\theta_2)\in\mathcal C.
		\]
		
		\item If $(G,\theta)\in\mathcal C$, then, for every $m\geq1$,
		\[
		\bigl(\mathbb Z_2^m\rtimes_\theta G,\widetilde\theta\bigr)
		\in\mathcal C.
		\]
	\end{enumerate}
	
	The members of $\mathcal C$ are called \emph{elementary type
		cyclotomic pro-$2$ pairs}.
\end{defn}




Motivated by Marshall's conjectures on Witt rings
\cite{marshall-book}, Efrat formulated the Elementary Type Conjecture
\cite{efrat-haran-94}. In terms of the class $\mathcal C$ defined
above, the version relevant to the present setting takes the following
form.


\begin{conjec}[Elementary Type Conjecture for maximal pro-$2$ Galois groups]
	Let $F$ be a field such that $G_F(2)$ is finitely
	generated as a pro-$2$ group. Then its canonical cyclotomic pro-$2$ pair $(G_F(2),\theta_F)$ belongs to $\mathcal C$.
\end{conjec}

Pythagorean fields with finitely many square classes and local fields
provide classical examples of fields $F$ for which
$(G_F(2),\theta_F)\in\mathcal C.$ 
We refer to \cite[Tables~5.1--5.3]{ware3} for the corresponding
cyclotomic pairs whose underlying pro-$2$ groups have at most three
generators.

\subsection{Cyclotomic semidirect products and valuations} 
\label{sec:semi-direct}

We now study fields whose maximal pro-$2$ Galois groups arise from
nontrivial cyclotomic semidirect products. Once the semidirect product
structure is fixed, only the underlying pro-$2$ groups and the
corresponding fields will be relevant, so we suppress the cyclotomic
characters from the notation. 

These structures are closely related to rigid elements and
$2$-henselian valuations, which we use to determine the Kaplansky
radical of the corresponding fields.


 


Let $F$ be a field. Recall that, for $a\in\dot F$, \[ D_F\langle1,a\rangle = \{x^2+ay^2\neq0\mid x,y\in F\}. \] 

\begin{defn}\label{def:rigid-birigid-elements} An element $a\in\dot F\setminus\bigl(\dot F^2\cup-\dot F^2\bigr)$ is called \emph{rigid} if \[D_F\langle1,a\rangle = \dot F^2\cup a\dot F^2.\] It is called \emph{birigid} if both $a$ and $-a$ are rigid. 
\end{defn} 

The existence of sufficiently many rigid elements gives rise to 
valuations; see \cite{arason1,engler3, ware2}. We first determine the
restriction imposed on the Kaplansky radical by a single rigid
element.

\begin{prop}\label{prop:fields-with-rigid} Let $F$ be a field containing a rigid element $a$. Then \[\bigl(R(F):\dot F^2\bigr)\leq2,\] and exactly one of the following alternatives occurs: 
	\begin{enumerate}[\upshape(a)]
		\item $R(F)=\dot F^2\neq\dot F.$ Thus, $F$ has trivial Kaplansky radical and $G_F(2)$ is not free pro-$2$. 
		\item The field $F$ is formally real and quasi-Pythagorean, and \[ R(F) = D_F\langle1,1\rangle = \dot F^2\cup a\dot F^2. \] 
		\item $R(F) = \dot F = \dot F^2\cup a\dot F^2$. In this case, $G_F(2)\cong\mathbb Z_2.$ 
\end{enumerate} 
\end{prop} 
\begin{proof}
Since $a$ is rigid,	
\begin{equation}\label{R(F)-rigid-indice-2} 	R(F)\subseteq D_F\langle1,a\rangle=\dot F^2\cup a\dot F^2.
\end{equation}
Therefore, $\bigl(R(F):\dot F^2\bigr)\leq2.$ If $a\notin R(F)$, it follows that $R(F)=\dot F^2.$ Moreover, since $R(F)\neq\dot F$, Proposition~\ref{prop:free-radical} implies that 
$G_F(2)$ is not free pro-$2$. This gives alternative~\textup{(a)}. 
Now assume that $a\in R(F)$. Then $\dot F^2\cup a\dot F^2\subseteq R(F)$. Together with \eqref{R(F)-rigid-indice-2}, it follows that 	\begin{equation}\label{eq:radical-rigid-two-classes}
R(F)=D_F\langle1,a\rangle=\dot F^2\cup a\dot F^2.
\end{equation}
Assume first that $R(F)=\dot F$. By Proposition~\ref{prop:free-radical}, $G_F(2)$ is free pro-$2$. From \eqref{eq:radical-rigid-two-classes} and Kummer theory, 	
\[\operatorname{rk}\bigl(G_F(2)\bigr) 	=\dim_{\mathbb F_2}\bigl(\dot F/\dot F^2\bigr)=1. 	\]
Hence $G_F(2)\cong\mathbb Z_2$. This proves~\textup{(c)}. It remains to consider the case $a\in R(F)$ and $R(F)\neq\dot F$. By
	\cite[Proposition~6.3(1), p.~451]{lam1},
	$D_F\langle1,a\rangle=D_F\langle1,1\rangle$. Again by 
	\eqref{eq:radical-rigid-two-classes}, $R(F) = D_F\langle1,1\rangle$, that is, $F$ is quasi-Pythagorean; see Remark~\ref{rem:chain}. Now observe that 	$-1\notin D_F\langle1,1\rangle$. Indeed, otherwise we have $D_F\langle1,1\rangle=\dot F$, contrary to
 $R(F)\neq\dot F$. It follows from 	\cite[Corollary~6.5(2), p.~452]{lam1} that $F$ is formally real. 	
This proves~\textup{(b)}.
\end{proof}

The following consequence will also be useful in the valuation-theoretic
arguments below.

\begin{cor}\label{cor:birigid-implies-trivia-radical}
Let $F$ be a field such that $(\dot F:\dot F^2)\geq4$ and $R(F)\neq\dot F.$ If $F$ contains a birigid element, then $R(F)=\dot F^2$.
\end{cor}
\begin{proof}
Let $a$ be birigid. By hypothesis, alternative~\textup{(c)} of
Proposition~\ref{prop:fields-with-rigid} cannot occur. 	Suppose first that $-1\in\dot F^2$. Then 	alternative~\textup{(b)} cannot occur either, since in that case $F$ 	is formally real. Hence alternative~\textup{(a)} holds, and therefore $R(F)=\dot F^2$. Suppose now that $-1\notin\dot F^2$. Then the square classes of $a$ and $-a$ are distinct. Since both $a$ and $-a$ are rigid, we have 
\[	R(F)\subseteq D_F\langle1,a\rangle\cap D_F\langle1,-a\rangle
	=	\bigl(\dot F^2\cup a\dot F^2\bigr)	\cap	\bigl(\dot F^2\cup(-a)\dot F^2\bigr)	=	\dot F^2.	\]
\end{proof}

The previous results show that rigid elements strongly restrict the
Kaplansky radical. We now turn to the valuation-theoretic counterpart
of this phenomenon. 

For a valued field $(F,v)$, let $\Gamma_v$ denote its value group,
$A_v$ its valuation ring, $m_v$ its maximal ideal, and
$k_v=A_v/m_v$ its residue field.

The next proposition describes how the position
of $R(F)$ relative to the group $(1+m_v)\dot F^2$ restricts the value
group and the residue field.

\begin{prop}{\cite[Proposition~2.4]{dario3}} \label{prop:radical-valued-field} Let $(F,v)$ be a valued field such that $R(F)\neq\dot F$ and $\operatorname{char}(k_v)\neq2.$ Then: 
	\begin{enumerate}[\upshape(a)]
		\item If $R(F)\nsubseteq(1+m_v)\dot F^2$, then $(\Gamma_v:2\Gamma_v)\leq2.$ 
		If $(\Gamma_v:2\Gamma_v)=2$, then $k_v$ is quadratically closed. 
		\item If $(1+m_v)\dot F^2\subsetneq R(F),$ then $\Gamma_v=2\Gamma_v.$ 
	\end{enumerate} 
\end{prop} 

To apply this proposition, we consider $2$-henselian valuations.

\begin{defn}\label{ex:2-henselian-fields} A valuation $v$ on the field $F$ is called \emph{$2$-henselian} if it has a unique extension to $F(2)$. If $\operatorname{char}(k_v)\neq2$, this is equivalent to \[ 1+m_v\subseteq\dot F^2; \] see \cite[Corollary~4.2.4]{engler2}.
\end{defn} 

Combining the preceding characterization with
Proposition~\ref{prop:radical-valued-field}, we obtain a criterion
for the triviality of the Kaplansky radical.

\begin{prop}\label{prop:radical-of-2hens} Let $(F,v)$ be a $2$-henselian valued field such that \[ \operatorname{char}(k_v)\neq2,\qquad \Gamma_v\neq2\Gamma_v, \qquad\text{and}\qquad R(F)\neq\dot F. \] Then $F$ has trivial Kaplansky radical. 
\end{prop} 
\begin{proof} Since $v$ is $2$-henselian, $(1+m_v)\dot F^2=\dot F^2.$ If $(\Gamma_v:2\Gamma_v)>2$, Proposition~ \ref{prop:radical-valued-field}\textup{(a)} implies $R(F)\subseteq\dot F^2.$ The same conclusion holds if $(\Gamma_v:2\Gamma_v)=2$ and $k_v$ is not quadratically closed. It remains to consider the case \[ (\Gamma_v:2\Gamma_v)=2 \qquad\text{and}\qquad \dot k_v=\dot k_v^2. \] Choose $x\in\dot F$ with $v(x)\notin2\Gamma_v$, and let $A_v^\times$ denote the group of units of the valuation ring. Then
	\[	\dot F=A_v^\times\dot F^2\cup xA_v^\times\dot F^2.	\]
	Since $\dot k_v=\dot k_v^2$ and $1+m_v\subseteq\dot F^2$, we have $A_v^\times\subseteq\dot F^2.$  Therefore $(\dot F:\dot F^2)=2.$ Since $R(F)\neq\dot F$, it follows that $R(F)=\dot F^2.$ 
\end{proof} 

Alternatively, for $(\dot F:\dot F^2)> 2,$ Proposition~\ref{prop:radical-of-2hens} follows from Corollary~\ref{cor:birigid-implies-trivia-radical}, since every element $a\in\dot F$ satisfying $v(a)\notin2\Gamma_v$ is birigid. 


The following theorem is a consequence of a deep result of Engler and
Nogueira \cite[Theorem~4.4]{engler4}, which establishes the existence
of a $2$-henselian valuation from the presence of a nontrivial normal
abelian subgroup in the maximal pro-$2$ Galois group.


\begin{thm}\label{eng-nog-thm} Let $F$ be a field such that
\[		G_F(2)\cong\mathbb Z_2^m\rtimes G_K(2),		\]
where the action defining the semidirect product is cyclotomic,
$m\geq1$, and $K$ is an intermediate field of $F(2)/F$ satisfying
$K\neq F(2)$. 
Assume that
$F$ is not a formally real Pythagorean field with 		
\[	G_F(2)\cong\mathbb Z_2\rtimes\mathbb Z/2\mathbb Z.		\]
Then $F$ admits a $2$-henselian valuation $v$ such that $\Gamma_v\neq2\Gamma_v$ and $\operatorname{char}(k_v)\neq2$.
\end{thm}
\begin{proof} Under the assumed semidirect product structure, $G_F(2)$ contains a nontrivial normal abelian subgroup isomorphic to $\mathbb Z_2^m$. The result therefore follows from \cite[Theorem~4.4]{engler4}.
\end{proof} 

\begin{cor}\label{cor:do-teorema-eng-nog}
	Let $F$ be a field such that
	\[
	G_F(2)\cong\mathbb Z_2^m\rtimes G_K(2),
	\]
where the action defining the semidirect product is cyclotomic, $m\geq1$, and $K$ is an
	intermediate field of $F(2)/F$ satisfying $K\neq F(2)$. Then
	$R(F)=\dot F^2$.
\end{cor}

\begin{proof}
	Suppose first that $F$ is a formally real Pythagorean field and that
	$G_F(2)\cong\mathbb Z_2\rtimes\mathbb Z/2\mathbb Z$. Since $F$ is
	Pythagorean, $D_F\langle1,1\rangle=\dot F^2$. As $R(F)\subseteq D_F\langle1,1\rangle$, we have that $R(F)$ is trivial. 
	Suppose now that the exceptional case does not occur. By
	Theorem~\ref{eng-nog-thm}, $F$ admits a $2$-henselian valuation $v$
	such that $\Gamma_v\neq2\Gamma_v$ and
	$\operatorname{char}(k_v)\neq2$. Since $G_F(2)$ is not free,  Proposition~\ref{prop:free-radical} gives 	$R(F)\neq\dot F$. The conclusion now follows from Proposition~\ref{prop:radical-of-2hens}.
\end{proof}

The preceding discussion also clarifies the relation between
$2$-henselian and pre-$2$-henselian fields. The latter, introduced in
\cite{dario3}, provide a natural source of fields with nontrivial
Kaplansky radical whose maximal pro-$2$ Galois groups admit free
pro-$2$ product decompositions.

\begin{ex}\label{ex:pre-2-henselian-fields} A valued field $(F,v)$ with 		$\operatorname{char}(k_v)\neq2$ is called 	\emph{pre-$2$-henselian} if 		\[1+m_v\subseteq R(F).\] Assume that $(F,v)$ is pre-$2$-henselian, $R(F)\neq\dot F$ and $\Gamma_v\neq2\Gamma_v.$  Since
$(1+m_v)\dot F^2\subseteq R(F),$ 
Proposition~\ref{prop:radical-valued-field}\textup{(b)} and the 	hypothesis $\Gamma_v\neq2\Gamma_v$ imply that the inclusion cannot be 
strict. Therefore, by \[R(F)=(1+m_v)\dot F^2.\] In particular, $R(F)=\dot F^2$ if and only if $v$ is $2$-henselian. 	Furthermore, \cite[Theorem~4.2]{dario3} shows that 
	\[	G_F(2)	\cong	G_L(2)\ast_2	\bigl(\mathbb Z_2^m\rtimes G_{k_v}(2)\bigr),	\]
	where $m=\dim_{\mathbb F_2}(\Gamma_v/2\Gamma_v),$ 	$L$ is an intermediate field of $F(2)/F$, $G_L(2)$ is a free pro-$2$ group, and $\dot L/\dot L^2	\cong 	R(F)/\dot F^2.$ Conversely, the existence of such a free pro-$2$ product
	decomposition characterizes pre-$2$-henselian fields.
\end{ex}

\subsection{The decomposition theorem in the elementary type case}
\label{sec:radical-ETG}
We now combine the results of the preceding subsections with the
behavior of the Kaplansky radical under free pro-$2$ products
established in Section~\ref{sec:R(F)-and-fre-products} to prove the
following result.



\begin{thm}\label{thm:claim2-ETG2}
Let $F$ be a field such that $(G_F(2),\theta_F)\in\mathcal C.$
Then there exist intermediate fields $L,H$ of $F(2)/F$ such that
\[	(G_F(2),\theta_F)	\cong	(G_L(2),\theta_L)\ast_2(G_H(2),\theta_H),
\]
where $G_L(2)$ is free pro-$2$ and $R(H)=\dot H^2.$
\end{thm}

\begin{proof}
	Let $\mathcal D$ be the subclass of $\mathcal C$ consisting of those
	pairs $(G,\theta)$ such that, for every field $E$ realizing
	$(G,\theta)$, there exist intermediate fields $L_E,H_E$ of
	$E(2)/E$ satisfying
	\[
	(G_E(2),\theta_E)
	\cong
	(G_{L_E}(2),\theta_{L_E})
	\ast_2
	(G_{H_E}(2),\theta_{H_E}),
	\]
	with $G_{L_E}(2)$ free pro-$2$ and $R(H_E)=\dot H_E^2.$
The basic pairs belong to $\mathcal D$. Indeed, if $G_E(2)$ is free
	pro-$2$, take $L_E=E$ and $H_E=E(2)$. If $E$ is a
	$2$-Demushkin field, take $L_E=E(2)$ and $H_E=E$ and apply
	Proposition~\ref{prop:demushkin-radical}. In both cases, the factor corresponding to $E(2)$ is the trivial
	cyclotomic pair $(\{1\},1)$. Thus the basic pairs belong to
	$\mathcal D$. We next show that $\mathcal D$ is closed under cyclotomic semidirect
products. Suppose that
\[
(G_E(2),\theta_E)
\cong
\bigl(\mathbb Z_2^m\rtimes_\theta G,\widetilde\theta\bigr),
\qquad (G,\theta)\in\mathcal D.
\]
If $G_E(2)$ is free pro-$2$, the first basic case applies. Otherwise,
we may assume that $G\neq\{1\}$. Indeed, if $G=\{1\}$, then
$G_E(2)\cong\mathbb Z_2^m$; since this group is not free pro-$2$, one
has $m\geq2$, and it may instead be written as
\[
\mathbb Z_2^m
\cong
\mathbb Z_2^{m-1}\rtimes\mathbb Z_2
\]
with trivial, and hence cyclotomic, action. Identify $G$ with its image in $G_E(2)$, and let $K$ be its fixed
field in $E(2)$. Then
\[
(G_K(2),\theta_K)\cong(G,\theta)
\qquad\text{and}\qquad
K\neq E(2).
\]
Hence Corollary~\ref{cor:do-teorema-eng-nog} gives
$R(E)=\dot E^2.$ Taking $L_E=E(2)$ and $H_E=E$ shows that the resulting pair belongs
to $\mathcal D$.
Finally, we prove that $\mathcal D$ is closed under free pro-$2$
products. Let $(G_1,\theta_1),(G_2,\theta_2)\in\mathcal D$ 
and suppose that
\[
(G_E(2),\theta_E)
\cong
(G_1,\theta_1)\ast_2(G_2,\theta_2)
\]
for some field $E$. Identify $G_1$ and $G_2$ with the corresponding
free factors of $G_E(2)$, and let $E_i$ be the fixed field of $G_i$
in $E(2)$, for $i=1,2$. Then
\[
(G_{E_i}(2),\theta_{E_i})
\cong
(G_i,\theta_i).
\]
Since $(G_i,\theta_i)\in\mathcal D$, there exist intermediate fields
$L_i,H_i$ of $E(2)/E_i$ such that
\[
(G_{E_i}(2),\theta_{E_i})
\cong
(G_{L_i}(2),\theta_{L_i})
\ast_2
(G_{H_i}(2),\theta_{H_i}),
\]
where $G_{L_i}(2)$ is free pro-$2$ and $R(H_i)=\dot H_i^2.$
 Regrouping the factors, we obtain
\[
\begin{aligned}
(G_E(2),\theta_E)
\cong{}&
\bigl(
(G_{L_1}(2),\theta_{L_1})
\ast_2
(G_{L_2}(2),\theta_{L_2})
\bigr)
\\
&\ast_2
\bigl(
(G_{H_1}(2),\theta_{H_1})
\ast_2
(G_{H_2}(2),\theta_{H_2})
\bigr).
\end{aligned}
\]
All the regroupings above are isomorphisms of cyclotomic pairs, by
the associativity and universal property of the free pro-$2$
product. Let $L$ and $H$ be the fixed fields in $E(2)$ of the two grouped
factors, respectively. Then
\[
(G_E(2),\theta_E)
\cong
(G_L(2),\theta_L)\ast_2(G_H(2),\theta_H),
\]
where
$G_L(2)\cong G_{L_1}(2)\ast_2G_{L_2}(2)$ and $G_H(2)
\cong
G_{H_1}(2)\ast_2G_{H_2}(2).$
Since a free pro-$2$ product of free pro-$2$ groups is free pro-$2$,
the group $G_L(2)$ is free pro-$2$. Moreover,
Theorem~\ref{radical-of-fre-product}\textup{(c)}, applied to
$G_H(2) \cong G_{H_1}(2)\ast_2G_{H_2}(2),$ gives $R(H)=\dot H^2.$
Therefore $\mathcal D$ is closed under free pro-$2$ products. 	Thus $\mathcal D$ contains the basic pairs and is closed under the  operations defining $\mathcal C$. By the minimality of $\mathcal C$, 	we have $\mathcal C=\mathcal D$, and the result follows.
\end{proof}

Combining Theorem~\ref{thm:claim2-ETG2} with
Theorem~\ref{thm:H90-decomposition}, we obtain the H-conjecture in
the elementary type case.

\begin{thm}\label{theorem90H} Let $F$ be a field such that $(G_F(2),\theta_F)\in\mathcal C.$ 	Let $K=F(\sqrt a),$ $a\in\dot F\setminus\dot F^2,$ and let $N=N_{K/F}:\dot K\longrightarrow\dot F$ 	be the norm map. Then 	\[	N^{-1}(R(F))=\dot F\,R(K).	\]
In particular, Conjecture~\ref{claimh90} holds for $F$.
\end{thm}
\begin{proof} By the definition of $\mathcal C$, the group $G_F(2)$ is finitely 	generated. Hence $\dot F/\dot F^2$ is finite, and therefore $F$ 	satisfies the finiteness hypothesis of Conjecture~\ref{claimh90}. 	By Theorem~\ref{thm:claim2-ETG2}, there exist intermediate fields $L,H$ of $F(2)/F$ such that \[	(G_F(2),\theta_F) 	\cong	(G_L(2),\theta_L)\ast_2(G_H(2),\theta_H),
	\] 	where $G_L(2)$ is free pro-$2$ and $	R(H)=\dot H^2.$
In particular, $G_F(2)\cong G_L(2)\ast_2G_H(2).$ Thus we just apply 	Theorem~\ref{thm:H90-decomposition}. 
\end{proof}



\section{Classical examples and the decomposition criterion}
\label{sec:revisiting}

We conclude by revisiting the classical constructions of fields with nontrivial Kaplansky radical due to Berman and Kula. Our purpose is twofold: first, to recall the numerical properties of the radicals arising in these constructions; second, to show that their maximal pro-$2$ Galois groups admit the free pro-$2$ product decomposition
required in Theorem~\ref{thm:H90-decomposition}.

Fields with nontrivial Kaplansky radical first arose in connection with the classification of Witt rings. Cordes \cite{cordes1} observed that such fields were needed to complete the classification of Witt rings of fields with eight square classes. Subsequently, Kula \cite{kula1} used quadratic form schemes to prove that, for arbitrary positive integers $n$ and $m$, there exists a field $F$ satisfying
\[(\dot F:R(F))=2^n \qquad\text{and}\qquad (R(F):\dot F^2)=2^m.\]
Berman \cite{berman} obtained further examples by studying nonreal
quadratic extensions of Pythagorean fields, while the question raised
by Cordes was later settled in \cite{kula2}. We begin by recalling
Berman's constructions.

	
\begin{ex}\cite[Theorem~2.3]{berman}\label{ex:berman1} Let $F$ be a formally real Pythagorean field such that the space of orderings $X(F)$ has cardinality $2^n$, where $n\ge2$, and 	$(\dot F:\dot F^2)=2^{n+1}.$ Let $K=F(\sqrt{-1}).$ Then \[	(\dot K:\dot K^2)=2^n	\qquad\text{and}\qquad 	(\dot K:R(K))=4.\] 
Consequently, $(R(K):\dot K^2)=2^{n-2},$ and hence $K$ has nontrivial Kaplansky radical when $n\geq 3$.
\end{ex}
	
Berman considerably generalized this construction, showing that the two quotients measuring the position of the Kaplansky radical can be controlled separately.
	
\begin{ex}\label{ex:berman2}\cite[Theorems 3.9 and 3.12]{berman} 	Let $F_1$ be a super-Pythagorean field with $2^n$ square classes and 	$2^{n-1}$ orderings, where $n\ge3$, and let $F_2$ be a SAP Pythagorean field with $2^m$ square classes and $m$ orderings, where $m\ge2$. Berman proved that there exists a field $F$ such that  
\[	\dot F/\dot F^{\,2}	\cong	\dot F_1/\dot F_1^{\,2}	\times	\dot F_2/\dot F_2^{\,2}.	\]
Setting $K=F(\sqrt{-1})$, he showed that \[\dot K/R(K)\cong\dot F_1/(\{\pm1\}\dot F_1^{\,2}) \qquad\text{and}\qquad R(K)/\dot K^{\,2}\cong\dot F_2/\dot F_2^{\,2}.\]
In particular,
\[
(\dot K:R(K))=2^{n-1}
\qquad\text{and}\qquad
(R(K):\dot K^2)=2^m.
\]
Thus both the codimension of $R(K)$ in $\dot K$ and the dimension of
$R(K)/\dot K^2$ can be prescribed within this family
\end{ex}


Kula \cite{kula1} constructed a family of fields for which the two
indices measuring the position of the Kaplansky radical can be
controlled separately. We recall the construction in the
Galois-theoretic form described in \cite[Section~4.2]{dario2}.

\begin{ex}\label{ex:kula}
Let $F_1$ and $L$ be fields with finitely many square classes. Assume that $R(F_1)=\dot F_1$ and that $L$ is not quadratically closed. Write $(\dot F_1:\dot F_1^2)=2^m,$ and, for an integer $n\geq1$, set $F_2=L((X_1))\cdots((X_n)).$ 
	Then there exists a field $F$ such that
$G_F(2)\cong G_{F_1}(2)\ast_2G_{F_2}(2).$	
Moreover,
\[	(R(F):\dot F^2)=2^m \quad \text{and} \quad  (\dot F:R(F))	=
	2^n(\dot L:\dot L^2). \]
\end{ex}

We conclude by showing that the classical constructions recalled above
satisfy the decomposition criterion of Theorem~\ref{thm:H90-decomposition}. 
Notice that the following result does not assert that the groups involved are of elementary type.


\begin{thm}\label{theo:examples-are-elementary}
Let $E$ be one of the fields described in Examples~\ref{ex:berman1}, \ref{ex:berman2}, or \ref{ex:kula}. Then
there exist intermediate fields $L,H$ of $E(2)/E$ such that
\[G_E(2)\cong G_L(2)\ast_2G_H(2),\]
where $G_L(2)$ is free pro-$2$ and $R(H)=\dot H^2.$ 
Consequently, the conclusion of
Theorem~\ref{thm:H90-decomposition} holds for $E$.
\end{thm}
\begin{proof} We consider the three constructions separately. First, consider Berman's construction in Example~\ref{ex:berman1}, where
\[	E=K=F(\sqrt{-1}).	\] 	Since $(\dot K:R(K))=4,$ every element of $\dot K\setminus R(K)$ is $R(K)$-rigid; see \cite[Lemma~2.1]{dario3}. Since $K$ is not formally real,
\cite[Theorem~2.2(1)]{dario3} implies that $K$ is pre-$2$-henselian. The required decomposition then follows from \cite[Theorem~4.2]{dario3}, as recalled in 		Example~\ref{ex:pre-2-henselian-fields}. 

Now consider Berman's construction in Example~\ref{ex:berman2}, where
\[E=K=F(\sqrt{-1}).	\] 	By \cite[Theorem~4.8]{dario3},
\[G_K(2)\cong G_{F_1(\sqrt{-1})}(2)	\ast_2 G_{F_2(\sqrt{-1})}(2)\ast_2\mathbb Z_2.\]
Since $F_2$ is SAP, $G_{F_1(\sqrt{-1})}(2)$ is a free pro-$2$ group by \cite[Remark~2]{ware1}. Furthermore, $F_2(\sqrt{-1})$ has trivial Kaplansky radical by
\cite{elmanlam2}. Regrouping the free factors, we obtain a decomposition
\[	G_K(2)\cong G_L(2)\ast_2G_H(2),	\]
where $G_L(2)$ is free pro-$2$ and $H$ has trivial Kaplansky radical. 

Finally, consider the Kula construction in
Example~\ref{ex:kula}. We have
\[
G_E(2)\cong G_{F_1}(2)\ast_2G_{F_2}(2),
\]
where $R(F_1)=\dot F_1$ and $F_2=L((X_1))\cdots((X_n)).$ By Proposition~\ref{prop:free-radical}, $G_{F_1}(2)$ is free
pro-$2$. The natural valuation on $F_2$ is $2$-henselian, has residue
field of characteristic different from $2$, and has a non-$2$-divisible
value group. Moreover, since $L$ is not quadratically closed,
$G_{F_2}(2)$ is not free pro-$2$. Hence
Proposition~\ref{prop:free-radical} gives $R(F_2)\neq\dot F_2.$ 
It now follows from Proposition~\ref{prop:radical-of-2hens} that $R(F_2)=\dot F_2^2.$ Thus, the displayed decomposition already has the required form.
\end{proof}

In the finite-dimensional setting, two further examples 
also satisfy the decomposition criterion. Pre-$2$-henselian fields
admit the required decomposition by
Example~\ref{ex:pre-2-henselian-fields}. Moreover, if $F$ is a
formally real quasi-Pythagorean field with $(\dot F:\dot F^2)<\infty,$ 
then Ware \cite[Corollary~1]{ware1} proved that $G_F(2)\cong\mathcal F\ast_2\mathcal H,$ where $\mathcal F$ is a free pro-$2$ group and $\mathcal H$ is
topologically generated by involutions. If $H$ denotes the fixed field
of $\mathcal H$ in $F(2)$, then $H$ is Pythagorean, and hence $R(H)=\dot H^2.$ 
Thus quasi-Pythagorean fields with finitely many square classes also
satisfy the hypotheses of Theorem~\ref{thm:H90-decomposition}.

Consequently, the decomposition criterion applies to all the
classical families with finitely many square classes considered in
this paper: pre-$2$-henselian fields, quasi-Pythagorean fields,  and
the constructions of Berman and Kula.

\backmatter

\bmhead{Acknowledgements}
I am deeply grateful to Antonio Jos\'e Engler for the great ideas for this paper and to J\'an Min\'a\v{c} for the endless mathematical discussions and all the inspiration for this work, partially written during my year at Western University. There, I also enjoyed working with Federico Pasini and Marina Palaisti.

\bibliography{sn-bibliography.bib}

\end{document}